\documentclass[11pt,dvips]{article}
 
 \usepackage{epsfig} 
 \usepackage{color,graphics}  
 \usepackage{latexsym,amsmath,amsfonts,amsthm,amssymb} 
 \usepackage{euscript,amsbsy} 
 \usepackage{graphicx}

\usepackage{latexsym}
\usepackage{amscd}
\usepackage{color,graphics}
\usepackage{graphicx}

\DeclareGraphicsExtensions{.eps}

\input colordvi    

\numberwithin{equation}{section}

\newcommand{\bE}{{\mathbb E}}
\newcommand{\bZ}{{\mathbb Z}}

\newcommand{\ga}{{\alpha}}

\newcommand{\gd}{{\delta}}
\newcommand{\gD}{{\Delta}}

\newcommand{\gl}{{\lambda}}
\newcommand{\gm}{{\mu}}

\newcommand{\gp}{{\pi}}

\newcommand{\gr}{{\varrho}}
\newcommand{\gs}{{\sigma}}

\newcommand{\gx}{{\xi}}

\newcommand{\gz}{{\zeta}}

\newtheorem{theo}{Theorem}[section]
\newtheorem{prop}[theo]{Proposition}
\newtheorem{lem}[theo]{Lemma}
\newtheorem{cor}[theo]{Corollary}

\newtheorem{defi}[theo]{Definition}
\newtheorem{rem}[theo]{Remark}
\newtheorem{ass}[theo]{Assumption}

\newcommand{\ba}{\begin{array}}
\newcommand{\ea}{\end{array}}
\newcommand{\bea}{\begin{eqnarray}}
\newcommand{\eea}{\end{eqnarray}}
\newcommand{\bead}{\begin{eqnarray*}}
\newcommand{\eead}{\end{eqnarray*}}
\newcommand{\be}{\begin{equation}}
\newcommand{\ee}{\end{equation}}
\newcommand{\bed}{\begin{displaymath}}
\newcommand{\eed}{\end{displaymath}}
\newcommand{\bl}{\begin{lem}}
\newcommand{\el}{\end{lem}}
\newcommand{\bp}{\begin{prop}}
\newcommand{\ep}{\end{prop}}
\newcommand{\bt}{\begin{theo}}
\newcommand{\et}{\end{theo}}

\newcommand{\bc}{\begin{cor}}
\newcommand{\ec}{\end{cor}}

\newcommand{\br}{\begin{rem}}
\newcommand{\er}{\end{rem}}
\newcommand{\bd}{\begin{defi}}
\newcommand{\ed}{\end{defi}}
\newcommand{\bass}{\begin{ass}}
\newcommand{\eass}{\end{ass}}

\begin{document}

\title{{   Queueing systems with pre-scheduled random arrivals}} 

\author{G. Guadagni$^{\dag}$ S. Ndreca$^{\ddag}$ B.Scoppola$^{\ddag}$}
\date{\today}
\maketitle

\begin{quote}
\centerline{\dag Department of Mathematics, University of Virginia}
\centerline{email: guadagni@virginia.edu}
\centerline{\ddag Dipartimento di Matematica, Universit\`a di Roma ``Tor
Vergata''} 
\centerline{e-mail:  ndreca@mat.uniroma2.it, scoppola@mat.uniroma2.it}
\end{quote}
\tableofcontents
\vspace{8mm}

\begin{abstract}
\noindent
We consider a point process $i+\xi_i$, where $i\in \bZ$ and the
$\xi_{i}$'s are i.i.d. random variables with variance $\sigma^{2}$.
This process, with a suitable rescaling of the distribution of
$\xi_i$'s, converges to the Poisson process in total variation for large
$\sigma$. We then study a simple queueing system
with our process as arrival process, and we provide a complete
analytical description of the system. 
Although the arrival process is very similar to the Poisson process,
due to negative autocorrelation the resulting queue is very
different from the Poisson case.
We found interesting
connections of this model with the statistical mechanics of
Fermi particles. 
This model is motivated by air traffic
systems. 
\end{abstract}

\bigskip
{\bf
Keywords:} Queueing system, air-traffic congestion, non Poissonian arrivals.
\bigskip


\renewcommand{\thefootnote}{\fnsymbol{footnote}}

\setcounter{footnote}{0}
\renewcommand{\thefootnote}{\arabic{footnote}}

\maketitle

\renewcommand{\thefootnote}{\fnsymbol{footnote}} \setcounter{footnote}{0} %
\renewcommand{\thefootnote}{\arabic{footnote}}

\section{Introduction}
The main aim of this paper is to define a stochastic point process 
to model the arrivals to a queueing system, and
to compare its features to the Poisson process. 

\noindent
It is well known that  the memoryless property
of the Poisson process simplifies many technical steps in the analysis
of queueing systems, but there are arrival processes where such an
assumption is not completely satisfied. In particular, we have in mind
air traffic models. In recent times the dramatic increase of air traffic stimulated a large
number of studies concerning the optimization of congestion management.
From the point of view of classical queueing theory the system is
difficult to study, mainly because it is hard even to define the
basic quantities of the theory. 
For instance it is clear that there is some congestion for  landing
aircrafts, since they have to follow some holding paths, 
but it is not easy to quantify the actual time spent in queue or even its 
instant length.
On the other hand, even assuming that the parameters of the system are known,
it is not clear what kind of point processes are suitable to  describe arrivals
and service times.
A common hypothesis in literature is to assume that arrivals are very well
modeled by a Poisson process. This assumption, to our
knowledge, goes back to the 70's when Dunlay and Horonjeff gave in
\cite{DH} a number of theoretical and statistical arguments 
to justify the Poissonian hypothesis,
and , since then, several other statistical studies
have supported the same results. Even recently, see \cite{tom}, a very careful
study of the interarrival times of aircrafts to major US airports shows
a small difference between the Poisson and the observed
distribution, i.e. the actual arrivals are slightly less random than
Poissonian ones, but the difference is quite small in all observed airports. 
On this ground, in various
papers, see for instance \cite{serra}, \cite{BEK} and \cite{BC} and reference therein, 
Poisson arrivals have
been assumed in the analysis of judicious management of service
times. It should be stressed that in all these papers  the statistical validation of the 
Poissonian hypothesis has been based
on computations on time scales  smaller than the intrinsic 
randomness of the system.

Stochastic models of aircraft arrivals based on statistical analysis
and on simulations have a long history. As a first attempt,
Barnett et al. \cite{barnett} studied  the arrivals to
Boston Logan Airport.
A version of the alternative model of arrivals we propose in this paper
was introduced and studied numerically in \cite{ball1}.
The model is refined in \cite{ball2}, where seasonal and daily
effects are taken into account to describe random
delays of departure times and,  with these
corrections, the model is quite accurate in its predictions.
The key feature of the model is a \textit{soft}  a-priori scheduling of arrivals: indeed, both in US and in
Europe, aircrafts 
are supposed to take off and to land by a schedule dictated by  
the capacity constraint of the runways, and by the assumption that each
aircraft would land in a very narrow time slot. However, on the day of operations,
an aircraft will be declared "on time" if it lands in a time interval larger than
ten times the original slot. In this sense the scheduling
should be considered "soft". 
The fact that arrivals are prescheduled clearly makes the Poissonian hypothesis
questionable, but this is usually neglected, on the basis of
the statistical studies mentioned above. However
the predictions of the queueing theory give in general very 
rough estimates of the actual queue length.
Moreover if we forecast a reduction of the
intrinsic variability of arrival times, which could be
achieved by various technical improvements
(e.g. a rescheduling closer to the actual arrival times, or an en-route control of the
paths of the aircrafts), we can not use Poissonian arrivals to describe
the system, because they depend only on  a single parameter $\rho$.

The process we study below is  an arrivals
model with two features. First, it shows a
pattern of arrivals very close to a Poisson process when we look at time
scales smaller than the standard deviation of aircraft delays,
second, it provides the distribution of arrivals  on time scales larger or
comparable to the standard deviation of aircraft delays.

Thus, the aim of this paper is an attempt to study more rigorously the features of
arrival process presented in \cite{ball1}, which we suitably generalize, 
and to understand its analytical properties.

Moreover, we show, both analytically and numerically, that the congestion related to
this process is very different from the congestion of a Poisson
process, on any time scale. This is due to the negative
autocorrelation of the process, as we prove explicitly.
It is worth to outline that the queueing models with Poisson arrivals
have in general probabilities to have $n$ customers in the queue
that decay much slower than the probabilities observed in the air traffic.
Our model gives a tail of the distribution much thinner, and more
similar to the observed data.

The analytical description of the system clarifies many interesting features of this kind of 
traffic: for heavy traffic the system has a long memory of the initial conditions; its description 
is obtained by the superposition of two processes, living on different time scales.
This give the possibility to investigate also systems with slowly variable traffic intensities.

The paper is organized as follows: in section 2 we describe our arrival process,
and we list some results on the comparison to the Poisson process. In section 3
we present a simplified computation, obtained neglecting the autocorrelation of the process.The congestion levels
according to this approximated process, assuming deterministic service (landing) times
and a single server (runway),
are quite different from the congestion according to Poisson arrivals. However
we show numerically that such approximation is bad for very congested systems, where
the actual level of congestion is not well described if the autocorrelation of the process is
neglected.
In section 4 we describe completely our queueing system 
at the price to enlarge suitably the state space of the Markov
chain describing it. It turns out that for our process
we have a finite value of the expected queue length even in the 
critical case
$\gr=1$, while
the Poisson queue diverges. Starting from the results 
on the critical case, we propose an
approximation scheme that works very well for highly congested
($\gr$ near to 1) systems. In this description a nice connection with
the statistical mechanics of Fermi gas emerges quite naturally.
Section 5 is devoted to conclusions and open problems.

\section{Description of the model: the arrival process}
In this section we want to introduce an arrival process, which we will call 
{\it pre-scheduled random arrivals} (PSRA) process, 
and to study its main features. The PSRA process is
defined as follows.
Let $\frac{1}{\lambda}$ be the expected interarrival time between two clients,
we define $t_i\in\mathbb{R}$ the actual arrival time of the $i$-th
client by
\be\label{2.1}
t_i={i\over\gl} +\xi_i \qquad i\in\mathbb{Z}
\ee
where $\xi_i$'s are i.i.d. random variables.

If the $\xi_{i}$'s are uniform, the model is the actual arrival times process introduced in
\cite{ball1} without cancellations and pop-ups. We will show later that
cancellations and pop-ups can be easily integrated into the process. 
From now on, we will assume that $\xi_i$'s have continuous
probability density $f^{(\gs)}_\xi(t)$ with variance $\sigma^{2}$, and we will set
without loss of generality $E(\xi_i)=0$,
since $E(\xi_i)\ne 0$ affects only the initial configuration of
the system. 
The main aim of this section is to compare the features of the PSRA
process to the Poisson process when $\sigma$ is large.
It is well known, e.g.  \cite[p.447]{Feller}, that the Poisson arrival process is defined by the fact
that probabilities $P_{j,j+1}(\gD t)=P(n(t+\gD t)=j+1|n(t)=j)$ of a "jump"
from the state $j$ to the state $j+1$ in the time interval $(t,t+\gD t]$ have the form  
\be
P_{j,j+1}(\gD t)=P^+(\gD t)=\gl\gD t+o(\gD t)
\ee
where $\gl$ is a constant independent of $t$ and
$j$; $\gl$ has the meaning of velocity of arrivals, i.e. denoting with $t_a$
the interarrival time, $E(t_a)=\frac{1}{\lambda}$.
For PSRA the probability 
$P(i,t,\gD t)$ that the $i$-th
client arrives in the time interval $(t,t+\gD t]$ is given by
\begin{align}
P(i,t,\gD t)&=P\left(t<{i\over\gl} +\xi_i<t+\gD t\right)=\\
&=P\left(t-{i\over\gl} <\xi_i<t+\gD t-{i\over\gl}\right)=
\int_{t-{i\over\gl}}^{t+\gD t-{i\over\gl}}f^{(\gs)}_\xi(x)dx
\end{align}
and, for small $\gD t$, it may be written as
\be\label{2.4}
P(i,t,\gD t)=f^{(\gs)}_\xi\left(t-{i\over\gl}\right)\gD t+o(\gD t)
\ee
By \eqref{2.4},  the probability
$P^+(t,\gD t)$ of a single PSRA arrival in the interval $(t,t+\gD t]$ is  
$$
P^+(t,\gD t)=\sum_{i\in\mathbb{Z}}P(i,t,\gD t)\prod_{j\ne i}(1-P(j,t,\gD t))=
$$
 $$=\sum_{i\in\mathbb{Z}}\left[f^{(\gs)}_\xi\left(t-{i\over\gl}\right)\gD t+o(\gD
t)\right] \exp\left(\sum_{j\ne i}\log\left[1-f^{(\gs)}_\xi\left(t-{j\over\gl}\right)\gD t+o(\gD
t)\right]\right) =$$
\be = \sum_{i\in\mathbb{Z}}\left[f^{(\gs)}_\xi\left(t-{i\over\gl}\right)\gD t+o(\gD
t)\right] \exp\left(-\sum_{j\ne i}\left[f^{(\gs)}_\xi\left(t-{j\over\gl}\right)\gD t+o(\gD
t)\right]\right) \ee

\begin{figure}[ht]   
\begin{center}
 \includegraphics[width=.80\textwidth,totalheight=2cm]{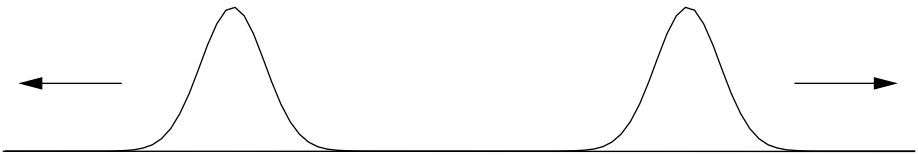}
\end{center}
\caption{A bimodal distribution with fixed shapes shifting to
infinity for $\gs\rightarrow \infty$.}\label{fig0}\nonumber
\end{figure}

\noindent
Hence up to the first order in $\gD t$ the rate of arrival $\gl(t)$ of the pre-scheduled random arrivals
is defined by 
\be \gl(t)=\sum_{i\in\mathbb{Z}}f^{(\gs)}_\xi\left(t-{i\over\gl}\right) 
\ee
This rate $\gl(t)$ is periodic in $t$ with period $\frac{1}{\lambda}$.
However we are interested in the dependence of $\gl(t)$ on $\gs$,
in particular when $\gs$ is large with respect to $\frac{1}{\lambda}$. 
To prove limit properties for our process, we have to specify the way we want to send $\sigma$ to infinity.
We will require the following \textit{scaling} property for the density $f^{(\sigma)}_{\xi} (t)$.
\bass
The probability density of $\xi$ has the form
\be\label{cond}
f^{(\gs)}_\xi(t,\gs^2)=\frac{1}{\sigma}f_\xi(t/\gs)
\ee 
i.e. it is the rescaling of a
well defined continuous density $f_\xi(t)$ with finite variance. We will also 
write $\max_{t\in{\mathbb {R}}} f_\xi(t)=M$.
\eass
This assumption is introduced in order to exclude pathological  
ways to send $\sigma$ to infinity, as, for instance, to 
consider a bimodal distribution with fixed maxima, see figure \ref{fig0}.

By our assumption, it follows that, in the limit $\gs$ very large the expression \be
R(\sigma,1/\gl):=\sum_{i\in\mathbb{Z}}{1\over\gl}f^{(\gs)}_\xi\left(t-{i\over\gl}\right) \ee
is the Riemann integral of the function $f^{(\gs)}_\xi(t)$.

\noindent 
For example, let $\xi$ be Gaussian $N (0,\sigma^{2} )$,
\[
R(\sigma,1/\gl)=
\sum_{i\in\mathbb{Z}}
{1\over\lambda }
\frac{1}{\sqrt{2\pi\sigma^{2}}}
e^{- \frac{(\lambda t-i)^{2}}{2\sigma^{2} \lambda^{2}} } =
\sum_{i\in\mathbb{Z}}
\frac{1}{\sqrt{2\pi}} 
e^{-\frac{1}{2}\left(\frac{t\lambda -i}{\lambda\sigma}
\right)^{2}}\frac{1}{\lambda \sigma} =
\sum_{i\in\mathbb{Z}}
\frac{1}{\sqrt{2\pi}} 
e^{-\frac{x_{i}^2}{2}}\Delta x 
\longrightarrow 1
\]
where $x_{i}=\frac{\lambda t-i}{\lambda\sigma}$ and $\Delta x=\frac{1}{\lambda \sigma }$
and the limit is for $\sigma \rightarrow \infty$.

\noindent
For any random variable rescaled in the above sense it is clear
that the result
\be
\lim_{\gs \to\infty}R(\sigma,1/\gl)=1
\ee
holds, and therefore, in the same limit,
\be\label{2.9}
\lim_{\gs\to\infty}\gl(t)=\lim_{\gs\to\infty}\lambda R(\sigma,1/\gl)= \gl
\ee
It is interesting,  for Gaussian $\xi$, to check numerically how fast
the limit is reached.
Table \ref{tab1} shows it. For 
simplicity, we set $\gl=1$.

\noindent
\begin{table}[ht]
\begin{center}
\begin{tabular}{|c|c|c|c|c|c|c|}
\hline
$\gs$&$\gl(0)$&$\gl(0.1)$&$\gl(0.2)$&$\gl(0.3)$&$\gl(0.4)$&$\gl(0.5)$\\
\hline
.2& 1.994726& 1.760407& 1.210523& 0.651951& 0.292114& 0.175283\\
.3& 1.340089& 1.274318& 1.103259& 0.894087& 0.726696& 0.663191\\
.4& 1.085005& 1.068767& 1.026261& 0.973729& 0.931237& 0.915008\\
.5& 1.014384& 1.011637& 1.004445& 0.995555& 0.988363& 0.985616\\
.6& 1.00164 & 1.001327& 1.000507& 0.999493& 0.998673& 0.99836\\
.7& 1.000126& 1.000102& 1.000039& 0.999961& 0.999898& 0.999874\\
.8& 1.000007& 1.000005& 1.000002& 0.999998& 0.999995& 0.999993\\
.9& 1.      & 1.      & 1.      & 1.      & 1.      & 1.\\
1.& 1.      & 1.      & 1.      & 1.      & 1.      & 1.\\
\hline
\end{tabular}

\begin{tabular}{|c|c|c|c|c|c|c|}
\hline
$\gs$&$\gl(0.6)$&$\gl(0.7)$&$\gl(0.8)$&$\gl(0.9)$&$\gl(1)$\\
\hline
.2& 0.292114& 0.651951& 1.210523& 1.760407& 1.994726\\
.3& 0.726696& 0.894087& 1.103259& 1.274318& 1.340089\\
.4& 0.931237& 0.973729& 1.026261& 1.068767& 1.085005\\
.5& 0.988363& 0.995555& 1.004445& 1.011637& 1.014384\\
.6& 0.998673& 0.999493& 1.000507& 1.001327& 1.00164\\
.7& 0.999898& 0.999961& 1.000039& 1.000102& 1.000126\\
.8& 0.999995& 0.999998& 1.000002& 1.000005& 1.000007\\
.9& 1.      & 1.      & 1.      & 1.      & 1.\\
1.& 1.      & 1.      & 1.      & 1.      & 1.\\
\hline

\end{tabular}
\caption{}\label{tab1}
\end{center}
\end{table}

\begin{figure}[ht]   
\begin{center}
 \includegraphics[width=.80\textwidth,totalheight=7cm]{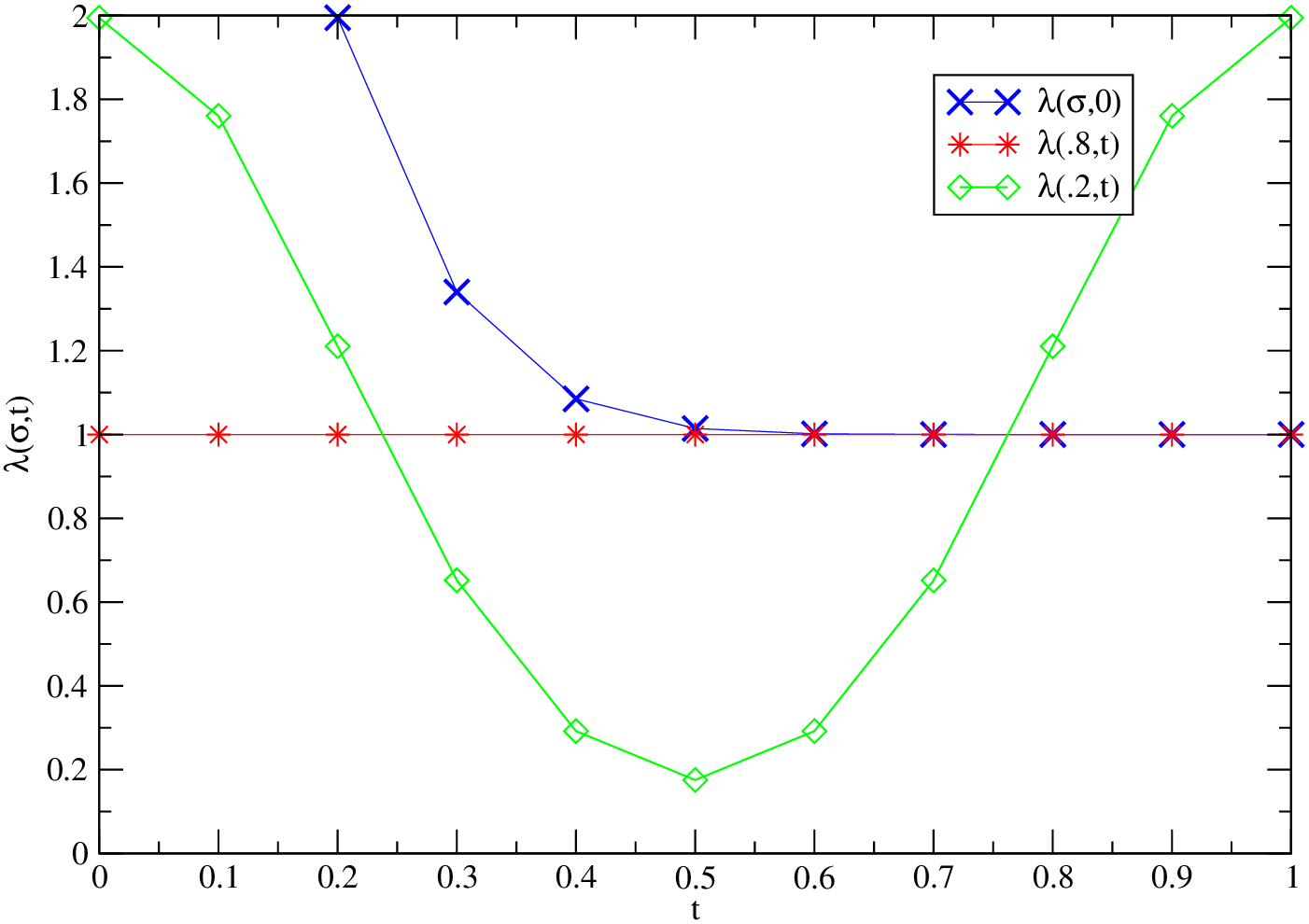}
\end{center}
\caption{Behavior of the function $\gl(\gs,t)$}\label{fig1}\nonumber
\end{figure}
The graph in figure \ref{fig1} shows that, in terms of rate
of arrivals, the pre-scheduled random arrivals approach the
Poisson process
when $\gs$ is suitably large. In particular for 
Gaussian variables with standard deviation $\gs$ of order  
$1/\gl$ or more we have that $\gl(t)$ is constant up to 6 digits.
Note that for applications mentioned in the introduction, we
do expect the standard deviation to be much larger than $1/\gl$.
Note also that the explicit structure of the density of $\xi$
does not play any particular role, and similar results may be obtained with
different distributions. 
\noindent
However it is clear that a small dependence
on $t$ is always present in the expression of $\gl(t)$, and hence it is
difficult to obtain a quantitative comparison between the pre-scheduled 
random arrivals and the Poisson process on this basis. Hence 
we look at the distribution of the random variable $n(t,t+T)$,
number of arrivals in the {\it finite} interval $(t,t+T]$.
Let us call $p_i(t,t+T)$ the probability that the $i$-th client arrives in the
interval $(t,t+T]$. Clearly
\be
p_i(t,t+T)=\int_t^{t+T}f^{(\sigma)}_\xi\left(x-{i\over\gl}\right)dx
\ee
Given the probabilities $p_i(t,t+T)$ we can write the generating function of the
random variable $n(t,t+T)$, and, defining $q_n^{(\gs)}=P(n(t,t+T)=n)$ 
we get 
\be\label{2.11}
q_n^{(\gs)}=\sum_{I=\{i_1,...,i_n\}}\prod_{i\in I}p_i(t,t+T)\prod_{j\notin I}(1-p_j(t,t+T))
\ee
where the sum runs over all the possible distinct subsets $I$ of indices
of cardinality $n$. By mean of this expression one obtains the generating function
\be\label{2.12}
q^{(\gs)}(z)=\sum_{n\ge 0}q_n^{(\gs)}z^n=
\prod_{i\in\mathbb{Z}}(1+(z-1)p_i(t,t+T))
\ee

\noindent
To take into account also the possibility of random independent {\it deletion} as in \cite{ball1},
let us outline here that a similar generating function can be introduced also
when each arrival has an independent probability $1-\gamma$ to be deleted, and
the complementary probability $\gamma$ to be an actual arrival. 
In other words, we construct the PSRA process for $i\in\mathbb{Z}$ and then for
each $i$ we cancel the corresponding $i$-th arrival with independent 
probability $1-\gamma$.
It is obvious that
in this case the generating function is
\be\label{2.12.1}
q^{(\gs)}_\gamma(z)=\sum_{n\ge 0}q_{\gamma,n}^{(\gs)}z^n=
\prod_{i\in\mathbb{Z}}(1+(z-1)\gamma p_i(t,t+T))
\ee

\noindent
The expressions \eqref{2.12}, \eqref{2.12.1} are exact, they give us all the information
on the distribution of $n(t,t+T)$, and they depend explicitly on $t$
and $T$.  However we want to study $q^{(\gs)}(z)$ and $q^{(\gs)}_\gamma(z)$
for large $\gs$, in the sense of the rescaling defined above, 
showing that they converge to a Poisson distribution with
parameter $\gl T$ and $\gamma\gl T$ respectively. The main idea is to exploit the fact that,
for large $\gs$, $p_i(t,t+T)$ goes to zero as ${1\over\gs}$. 

\noindent
We now prove the following results.
\bl \label{max}
\be
\max_i p_i(t,t+T)\le\frac{{const(T)}}{\gs}
\ee
\el
\begin{proof}
\be
p_i(t,t+T)=\int_t^{t+T}f^{(\gs)}_\xi\left(x-{i\over\gl}\right)dx=
\int_{t-\frac{i}{\gl}}^{t-\frac{i}{\gl}+T}f^{(\gs)}_\xi(s)ds={1\over\gs}
\int_{{t}-\frac{i}{\gl}}^{{t}-\frac{i}{\gl}+{T}}
f_\xi\left({s\over\gs}\right)ds
\ee
by the Intermediate Value Theorem 
\be
p_i(t,t+T)={1\over\gs}f_\xi\left({s_i^{*}\over\gs}\right){T}\le{MT\over\gs} 
\ee
where
\[
\frac{s_i^{*}}{\sigma } \in \left(t-\frac{i}{
\lambda},t-\frac{i}{
\lambda}+T \right)
\]
\end{proof}
\noindent 
Now we will use lemma \ref{max} to bound the generating function
\begin{gather}
q^{(\sigma)}(z)=\exp\left[\sum_{i\in\mathbb{Z}}\ln(1+(z-1)p_i(t,t+T))\right]= \\
=\exp\left[(z-1)\sum_{i\in\mathbb{Z}}p_i(t,t+T)\left(1+(z-1)p_i(t,t+T)\int_0^1ds
\frac{s}{(1+(z-1)(1-s)p_i(t,t+T))^2}\right)\right] \label{gen}
\end{gather}

\bl \label{intensita'}
With $p_i(t,t+T)$ defined as above, the sum in \eqref{gen} converges
to $\lambda T$
\be
\lim_{\gs\to\infty}\sum_{i\in\mathbb{Z}}p_i(t,t+T)
\left(1+(z-1)p_i(t,t+T)\int_0^1ds
\frac{s}{(1+(z-1)(1-s)p_i(t,t+T))^2}\right)=\gl T
\ee
\el
\begin{proof}
First we prove that 
\be\label{leading}
\lim_{\gs\to\infty}\sum_{i\in\mathbb{Z}}p_i(t,t+T)=\gl T.
\ee
Let us define
$T:=\frac{K+\gD T}{\gl}$, where  $K\in\mathbb{Z^+}$ and 
$0\le\gD T<1$. Then we can write
\[
\sum_{i\in\mathbb{Z}}p_i(t,t+T)=
\sum_{i\in\mathbb{Z}}\int_{t-\frac{i}{\gl}}^{t-\frac{i}{\gl}+T}f^{(\gs)}_\xi(s)ds=
\sum_{i\in\mathbb{Z}}\int_{t-\frac{i}{\gl}}^{t+\frac{K-i}{\gl}+
{\gD T\over \gl}}f^{(\gs)}_\xi(s)ds=
\]
\be \label{intero}
=
\sum_{i\in\mathbb{Z}}\int_{t-\frac{i}{\gl}}^{t+\frac{K-i}{\gl}}f^{(\gs)}_\xi(s)ds+
\sum_{i\in\mathbb{Z}}\int_{t+\frac{K-i}{\gl}}^{t+\frac{K-i}{\gl}+
{\gD T\over \gl}}f^{(\gs)}_\xi(s)ds
\ee
The first term on the right hand side of \eqref{intero} is $K$. 
Let  $i=mK+l$, where $l\in\mathbb{Z^+}$ and $m\in\mathbb{Z}$,
\be
\sum_{i\in\mathbb{Z}}\int_{t-\frac{i}{\gl}}^{t+\frac{K-i}{\gl}}f^{(\gs)}_\xi(s)ds=
\sum_{l=0}^{K-1}\sum_{m\in\mathbb{Z}}
\int_{t-\frac{mK+l}{\gl}}^{t-\frac{(m-1)K+l}{\gl}}f^{(\gs)}_\xi(s)ds=
\sum_{l=0}^{K-1}\int_{\mathbb{R}}f^{(\gs)}_\xi(s)ds=K
\ee
The second term on the right hand side of \eqref{intero} 
converges to $\gD T$ for $\gs\rightarrow\infty $:
\be
\sum_{i\in\mathbb{Z}}\int_{t+\frac{K-i}{\gl}}^{t+\frac{K-i}{\gl}+
{\gD T\over \gl}}f^{(\gs)}_\xi(s)ds=
\sum_{i\in\mathbb{Z}}\int_{t+\frac{i}{\gl}}^{t+\frac{i}{\gl}+
{\gD T\over \gl}}f^{(\gs)}_\xi(s)ds=
\sum_{i\in\mathbb{Z}}{1\over\gs}\int_{{t}+
\frac{i}{\gl}}^{{t}+\frac{i}{\gl}+{\gD T\over \gl}}
f_\xi\left({s\over\gs}\right)ds
\ee
and, by the Intermediate Value Theorem we get
\be
\sum_{i\in\mathbb{Z}}{1\over\gs}\int_{t+\frac{K-i}{\gl}}^{t+\frac{K-i}{\gl}+
{\gD T\over \gl}}f_\xi\left({s\over\gs}\right)ds=
\sum_{i\in\mathbb{Z}}{1\over\gs}f_\xi\left({s_i^{*}\over\gs}\right){\gD
T\over\gl}=
\gD T\sum_{i\in\mathbb{Z}}f_\xi\left({s_i^{*}\over\gs}\right){1\over\gl\gs}
\longrightarrow \gD T 
\ee
as $\gs\to\infty$, where the sum on the last equality 
is the Riemann sum of $f_\xi(t)$. This ends the proof of \eqref{leading}.
In order to complete the lemma we need to show that, uniformly in $i$,
$$
\lim_{\gs\to\infty}(z-1)p_i(t,t+T)\int_0^1ds
\frac{s}{(1+(z-1)(1-s)p_i(t,t+T))^2}=0
$$
but this follows from lemma \ref{max} and from the fact that
$$
(z-1)\int_0^1ds\frac{s}{(1+(z-1)(1-s)p_i(t,t+T))^2}\le C
$$
for any $p_i(t,t+T)<1/2$ and $|z|\leq1$.

\end{proof}

\bl \label{generating}
Let $q(z)=\exp({\gl T(z-1)} )$ be
the probability generating function of the Poisson 
random variable $\gz$ with intensity $\gl T$, and 
$q_\gamma(z)=\exp({\gamma\gl T(z-1)} )$ be
the probability generating function of the Poisson 
random variable $\gz$ with intensity $\gamma\gl T$, then  
\be \label{continuita}
\lim_{\gs\to\infty}q^{(\gs)}(z)=q(z);\quad \lim_{\gs\to\infty}q_\gamma^{(\gs)}(z)=q_\gamma(z)
\ee
\el
\begin{proof} 
Follows immediately  from lemma \ref{intensita'}.
\end{proof}
\bt
If $q^{(\gs)}(z)\longrightarrow q(z)$, then 
$\sum_{n=0}^{\infty}|q_n^{(\gs)}-q_n|\longrightarrow 0$ as $\gs\to\infty$.
The same result holds for the arrivals with random deletions.
\et
\begin{proof}
The proof follows from the continuity theorem for probability
generating function see Feller \cite[p.280]{Feller}.
\end{proof}
\vglue.5truecm
\noindent
Hence the PSRA process converges in distribution
to the Poisson process in total variation norm, and the same is true
for PSRA process with independent random deletions.

\noindent
In order to show that the process has negative autocorrelation, we
will compute the 
expected value, the variance $Var(n)$ of the number $n$ of arrivals in a 
time slot $(t,t+T]$, and the covariance $Cov (n_{1},n_{2})$,  
where $n_{1}$ and $n_{2}$ are the numbers of arrivals in $(t,t+T]$
and $(t+T,t+2T]$, respectively. We present the explicit computations in the case
of simple PSRA process, but the same results are true with obvious modifications
for PSRA process with independent random deletions.\\
Let $\chi_{i} (t_{i}\in (t,t+T])$ be the characteristic function of
the event ``client $i$ arrives in the interval $(t,t+T]$'', so
that $\bE (\chi_{i})=p_{i}(t,t+T)$, then the expected number of arrivals
in a time slot $(t,t+T]$ is
\[
\bE (n)=\bE \left(\sum_i \chi_{i} \right) = \sum_{i}\bE
(\chi_{i})=\sum_{i} p_{i}(t,t+T) 
\]
and also
\begin{align*}
\bE (n^{2})&= \bE \left(\sum_i \chi_{i} \sum_j \chi_{j}  \right)=
\bE \left(\sum_i \chi_{i} + \sum_{i\neq j} \chi_{i}\chi_{j}
\right) =\\
&= \sum_{i} p_{i}(t,t+T) + \sum_{i\neq j} p_{i}(t,t+T)p_{j}(t+T,t+2T) \\
&= \sum_{i} p_{i}(t,t+T) +
\left( \sum_{i} p_{i}(t,t+T) \right)^{2} - \sum_{i} (p_{i}(t,t+T))^{2}
\end{align*}
Then the variance 
is:  
\[
Var(n)= \bE (n^{2})- (\bE (n))^{2} = \sum_{i} p_{i}(t,t+T) -
\sum_{i} (p_{i}(t,t+T))^{2} = \sum_{i} p_{i}(t,t+T) (1-p_{i}(t,t+T))
\]
and we see again that $Var(n) \rightarrow \lambda T$ in the limit
$\gs\rightarrow\infty$.
Finally, 
let us define $\chi_{i}^{(1)}:=\chi_{i} (t_{i}\in (t,t+T])$
and $\chi_{i}^{(2)}:=\chi_{i} (t_{i}\in (t+T,t+2T])$ 
\begin{align*}
\bE (n_{1}n_{2})&= \bE
\left(\sum_{i}\chi^{(1)}_{i}\sum_{j}\chi^{(2)}_{j} \right) =
\bE ( \sum_{i\neq j} \chi_{i}^{(1)}\chi_{j}^{(2)}) = 
\sum_{i\neq j} \bE ( \chi_{i}^{(1)})\bE ( \chi_{j}^{(2)}) =\\
&= \sum_{i\neq j} p_{i}(t,t+T)p_{j}(t+T,t+2T) \\
&= \sum_{i,j}
p_{i}(t,t+T)p_{j}(t+T,t+2T) - \sum_{i} p_{i}(t,t+T)p_{i}(t+T,t+2T)
\end{align*}
so that 
\[
Cov (n_{1},n_{2})=\bE (n_{1}n_{2}) - \bE (n_{1})\bE (n_{2}) = 
-\sum_{i} p_{i}(t,t+T)p_{i}(t+T,t+2T)
\]
A negative covariance means that $n_{1}$ and $n_{2}$ are
inversely correlated, as we should expect in our arrival
model: a congested time slot should be followed or preceded by a
slot with lower than expected arrivals. 
Moreover, this is a clear indication that the hypothesis of
independence for 
$n_{1}$ and $n_{2}$, numbers of arrivals in different time slots,
is not correct, unless we are in the limit $\sigma \rightarrow \infty$.

\section{Queueing systems with PSRA process: independence approximation}

In this section we want to try to use the classical results of queueing
theory for a system in which the arrivals are described in terms of our PSRA,
there is a single server and the service time is deterministic.
For the air traffic applications the deterministic service (landing) times
are obviously an approximation, but neglecting the mix of aircrafts 
the actual landing times have a low variability.

\begin{figure}[t]\label{intens}   
\begin{center}
 \includegraphics[width=.80\textwidth,totalheight=6.0cm]{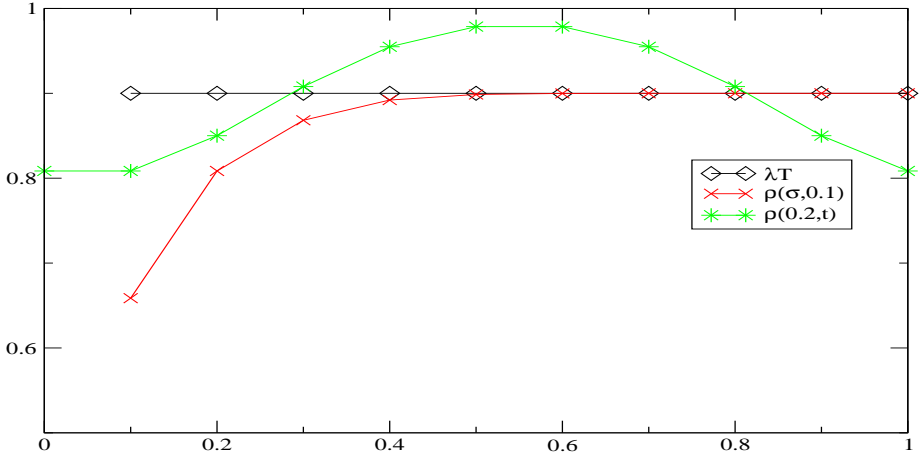}
\end{center}
\caption{Behavior of the function $\gr(\gs,t)$.
On the $x$ axis we have time $t$ for $\gr(0.2,t)$ and standard deviation $\gs$
for $\gr(\gs,0.1)$.}\label{fig2}
\nonumber
\end{figure}

In order to study our queueing 
process we set a service time $T$ and we define the instant traffic intensity
$\gr(\gs,t)=E(n(t,t+T))$. In fig. \ref{fig2} and table \ref{tab3} we report
numerical results for the convergence  of $\gr(\gs,t)$ to $\gl T$, granted by lemma \ref{intensita'}.
For simplicity we consider $\gx$ Gaussian, and $\gl=1$.  In this case $\gr(\gs,t)$ converges as soon as $\gs$
gets close to 1.
\begin{table}
\begin{center}
\begin{tabular}{|c|c|c|c|c|c|c|}
\hline
$\gs$&$T$&$\gr(\gs,0)$&$\gr(\gs,0.1)$&$\gr(\gs,0.2)$&$\gr(\gs,0.3))$&$\gr(\gs,0.4))$\\
\hline
.2&.9&0.808534&0.808534& 0.850089& 0.907951& 0.954826\\
.3&.9&0.868214&0.868214& 0.88048& 0.900153& 0.919615\\
.4&.9&0.892048&0.892048& 0.895086& 0.900001& 0.904914\\
.5&.9&0.898654&0.898654& 0.899168& 0.9& 0.900832\\
.6&.9&0.899847&0.899847& 0.899905& 0.9& 0.900095\\
.7&.9&0.899988&0.899988& 0.899993& 0.9& 0.900007\\
.8&.9&0.899999&0.899999& 0.9& 0.9& 0.9\\
.9&.9&0.9         &0.9 & 0.9& 0.9& 0.9\\
1.& .9&0.9      &0.9   & 0.9& 0.9& 0.9\\
\hline
$\gs$&$T$&$\gr(\gs,0.5)$&$\gr(\gs,0.6)$&$\gr(\gs,0.7)$&$\gr(\gs,0.8)$&$\gr(\gs,0.9)$\\
\hline
.2&.9&0.9786 &0.9786 &0.954826 &0.907951 &0.850089 \\
.3&.9&0.931537 &0.931537 &0.919615 &0.900153 &0.88048 \\
.4&.9&0.907951 &0.907951 &0.904914 &0.900001 &0.895086 \\
.5&.9&0.901346 &0.901346 &0.900832 &0.9 &0.899168 \\
.6&.9&0.900153 &0.900153 &0.900095 &0.9 &0.899905 \\
.7&.9&0.900012 &0.900012 &0.900007 &0.9 &0.899993 \\
.8&.9&0.900001 &0.900001 &0.9 &0.9 &0.9 \\
.9&.9&0.9 &0.9 &0.9 &0.9 &0.9 \\
1.& .9&0.9 &0.9 &0.9 &0.9 &0.9 \\
\hline
\end{tabular}
\caption{}
\label{tab3}
\end{center}
\end{table}
\noindent

We want to compare
the average queue size in $M/D/1$ queueing system (Poisson arrivals)
with a queueing system in which the arrivals are described
in terms of  PSRA.
To do this we have to recall some standard results in queueing theory.
Assuming to have a probability $Q_n$  to have $n$ arrivals in a 
service time slot, and assuming the variables $n$ to be i.i.d, 
our system is described by the so-called discrete time $GI/D/1$ queueing model.

It is well known, see e.g.\cite{Tijms}, that the stationary probabilities for
the discrete time $GI/D/1$ queueing model are given by
\be \label{stationary}
\begin{split}
&P_0=(P_0+P_1)Q_0\\
&\vdots \\
&P_n=P_0Q_n+\sum_{k=1}^{n+1}P_k Q_{n-k+1}
\end{split}
\ee

The corresponding generating function is 
\be\label{genf}
P(z)=\frac{P_0(1-z)}{1-{z\over Q(z)}} 
\ee
\\*
In the case of Poisson arrivals
with traffic intensity $\gr$,
$Q(z)=q(z)= \exp(\gr(z-1))$. Denoting by $N$ the
average queue size, 
after straightforward computations we get 
\be \label{poiss}
N=\frac{\gr(2-\gr)}{2(1-\gr)}
\ee
Consider now the PSRA process. In this case we can try to compute
\eqref{genf} by means of the generating function \eqref{2.12}. This is obviously an
approximation, since for PSRA arrivals, as it has been shown
in Section 2, the number of arrivals in subsequent time slots are
not independent.

However, neglecting the autocorrelation, we have
that $Q(z)=q^{(\gs)}(z)$, and denoting by $N(\gs,t)$
the average queue size we find

\be \label{valore}
N(\gs,t)=\frac{2\sum_{i\in\mathbb{Z}}p_i(t,t+T)-
(\sum_{i\in\mathbb{Z}}p_i(t,t+T))^2-
\sum_{i\in\mathbb{Z}}p_i^2(t,t+T)}
{2(1-\sum_{i\in\mathbb{Z}}p_i(t,t+T))}
\ee

For $\sigma$ large $N{(\gs,t)}$ becomes independent
of $t$, and it converges to $N$ by \eqref{leading}. 
Table \ref{tab4} shows that for Gaussian $\xi$ and $\lambda=1$ the convergence is quite fast.

\begin{table}[h]
\begin{center}
\begin{tabular}{|c|c|c|c|c|c|c|c|}
\hline
$\gs$&$T$&$N(\gs,0)$&$N(\gs,0.1)$&$N(\gs,0.2)$&$N(\gs,0.3)$&$N(\gs,0.4)$&$N(\gs,0.5)$\\
\hline
.1&.9&0.89105&0.89105&1.00493&1.04024&1.02267&1.00905\\
.2&.9&1.61425&1.61425&1.58187&1.51872&1.42902&1.32201\\
.3&.9&2.26812&2.26812&2.21399&2.10656&1.95949&1.83453\\
.4&.9&2.75253&2.75253&2.68673&2.57205&2.44587&2.36133\\
.5&.9&3.03548&3.03548&2.9955 &2.92993&2.86327&2.82151\\
.6&.9&3.24502&3.24502&3.23019&3.20614&3.18205&3.16714\\
.7&.9&3.43207&3.43207&3.42809&3.42165&3.41521&3.41123\\
.8&.9&3.59488&3.59488&3.59405&3.5927 &3.59134&3.59051\\
.9&.9&3.73131&3.73131&3.73117&3.73094&3.73071&3.73056\\
1.&.9&3.84462&3.84462&3.8446 &3.84457&3.84454&3.84452\\
\hline
\end{tabular}
\begin{tabular}{|c|c|c|c|c|c|c|}
\hline
$\gs$&$T$&$N(\gs,0.6)$&$N(\gs,0.7)$&$N(\gs,0.8)$&$N(\gs,0.9)$&$N(\gs,1)$\\
\hline
.1&.9&1.00905&1.02267&1.04024&1.00493&0.89105\\
.2&.9&1.32201&1.42902&1.51872&1.58187&1.61425\\
.3&.9&1.83453&1.95949&2.10656&2.21399&2.26812\\
.4&.9&2.36133&2.44587&2.57205&2.68673&2.75253\\
.5&.9&2.82151&2.86327&2.92993&2.9955 &3.03548\\
.6&.9&3.16714&3.18205&3.20614&3.23019&3.24502\\
.7&.9&3.41123&3.41521&3.42165&3.42809&3.43207\\
.8&.9&3.59051&3.59134&3.5927 &3.59405&3.59488\\
.9&.9&3.73056&3.73071&3.73094&3.73117&3.73131\\
1.&.9&3.84452&3.84454&3.84457&3.8446 &3.84462\\
\hline
\end{tabular}
\caption{}
\label{tab4}
\end{center}
\end{table}
\noindent

\begin{figure}[ht] 
\begin{center}
\includegraphics[width=.80\textwidth,totalheight=6cm]{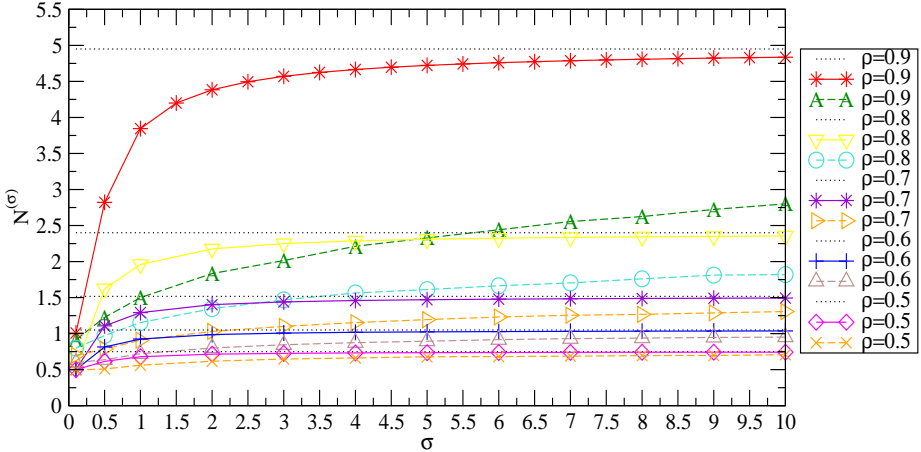}
\end{center}
\caption{Behavior of the function $N{(\gs,0.5)}$, for different values of $\gr$.
Dotted lines refer to Poisson arrivals, continuous lines refer to approximation
\eqref{valore}, dashed lines refer to simulations. The simulations are run for a time sufficiently
long to have fluctuations on the result negligible in the scale of the figure.}
\label{aver1}
\end{figure}

The results obtained by the formulas above are an approximation, because
we neglected the (negative) autocorrelations, and we have to see
when this approximation is reliable. As a matter of fact the PSRA process is easy to
implement for numerical simulation; hence we can
compare  the PSRA average queue size $N(\sigma,t)$ obtained by numerical simulations to \eqref{valore} and \eqref{poiss}.  
In figure~\ref{aver1} $N(\sigma,t)$ is plotted as a function of $\gs$, for
different values of $\gr=0.5,0.7,0.9$, and $t=0.5$. 
The dotted straight lines represent $N$ obtained by   
\eqref{poiss} for different values of $\gr$. 
As we can see from the graph, values of $N(\sigma,t=0.5)$ for fixed 
$\gr$ given by \eqref{valore} are larger than the corresponding ones obtained by simulation. 
Moreover,  this
overestimate becomes very important when $\gr$ increases. Hence, as it was easy to guess,
the negative autocorrelation plays an important role in the system when the traffic intensity
becomes large. For air traffic applications $\gr$ near to the critical value $\gr=1$
is the interesting case. 
\noindent

\section{Queueing systems with PSRA process: autocorrelated arrivals}

As it is clear from the results of the previous section, neglecting the autocorrelation
the computed average queue length is grossly overestimated
in the interesting cases. If we want to describe the system only by the length of the queue,
the presence of autocorrelation implies the loss of Markov property.
In this section we show that if we enlarge suitably the state space we may keep the Markov property, and describe completely the autocorrelation.
With this description some
interesting features of the system are clarified, but at
the moment we are able to compute explicitly the quantities of interest with some
approximations. Such approximations, however, turn out to give almost negligible errors.

\noindent 
To simplify the analytical treatment of the system, we will consider from now on
densities $f^{(\gs)}_\xi(t)$ of the random i.i.d. variables $\xi_i$ that are compact support,
i.e. such that $f^{(\gs)}_\xi(t)=0$ for $|t|>L$ for some $L<\infty$. We are setting $\gl=1$,
and we take $L\in\mathbb{N}$.
This implies that at a certain discrete time $j$ the $i$'th customer is certainly arrived to the
system for all $i\le j-L$, while for all $i\ge j+L$ it is certainly not yet arrived.
Hence to completely describe the state of the system we have to specify, beside the
number $n$ of customers waiting in queue right before the service at time $j$
is delivered, also a {\it finite} set $I_j$ of $i$'s,  $I_j\subset\{j-L+1,...,j+L-1\}$, that are the customers
that are already arrived at the service at time $j$. Note that the customers in the set $I_j$
are not necessarily already served at time $j$, or, in other words, the set $I_j$ is the set of 
the customers with indices in $\{j-L+1,...,j+L-1\}$ that are in the queue at time $j$, or that
are already served at time $j$. Note also that $0\le|I_j|\le2L-1$. Finally, we want to outline
that due to the independence of the $\xi$'s $I_{j+i}$ is independent of $I_j$
for all $i\ge 2L$.

\noindent
We will treat first the case $\gr=1$, or in other words, the case
$\gl=T=1$ in \eqref{2.12}. This special case is important
for several reasons. First, we will prove that for PSRA arrivals
the system has a {\it finite} average queue length, showing that,
even if the PSRA process tends in distribution to the Poisson process,
for finite variance of the $\xi$'s the two systems are deeply different.
Second, we will show that in the $\gr=1$ case  there is a conserved quantity
in the system, when the stationary distribution is reached. Third, it is
possible, using an interest interpretation of the system in terms of
Fermi statistics, to compute the (very long) time needed to the system to reach
the stationary distribution. Fourth, and maybe more important, on the
basis of this computation it is possible to approximate
efficiently the distribution of the length of the queue even for $\gr<1$.

\noindent
Hence, we fix $\gr=1$ and we start from the obvious relation
\be \label{5.1}
n(j+1)=n(j)-(1-\gd_{n(j)0})+m(j)
\ee
where $n(j)$ is the length of the queue immediately before the service at time $j$,
$m(j)$ is the number of customers arrived in the time slot $[j,j+1)$, and the term
$(1-\gd_{n(j)0})$ indicates the fact that if there is some customer in the queue at
time $j$, i.e. $n(j)>0$, the first of the queue is served, while if $n(j)=0$ then 
$n(j+1)=m(j)$.

\noindent
Now we observe that with our notations we can write
\be \label{5.2}
m(j)= |I_{j+1}|-|I_j|+1
\ee
This relation
can be shown as follows: the total number $na(j)$ of
customers arrived to the service from a certain fixed time, say from time $1$, to time $j$, 
is obviously $na(j)= j-L+|I_j|$,
because all the customers $k$ up to customer $j-L$'th are already arrived, due to the 
compactness of the support of $f^{(\gs)}_\xi(t)$, while for $k>j-L$ the number of arrived customers 
is $|I_j|$ by definition. Hence $m(j)=na(j+1)-na(j)= j+1-L+|I_{j+1}|-j+L-|I_j|= |I_{j+1}|-|I_j|+1$.
Putting \eqref{5.2} into\eqref{5.1} we obtain
\be \label{5.3}
n(j+1)=n(j)+|I_{j+1}|-|I_j|+\gd_{n(j)0}
\ee
This relation shows that the quantity $\ga(j)=n(j)-|I_j|$ is constant during a busy period,
and it increases by 1 at the end of each busy period. This implies that the stationary distribution
is reached once $\ga>0$. If the initial value of $\ga$ is strictly positive, the value $n(j)=0$
is never realized, and then $\ga$ remains constant and
\be \label{5.4.1}N=E(n)=\ga+E(|I|)
\ee
If the initial value of $\ga$ is 0 or it is
negative, a sequence of busy periods is realized, giving in the end the value $\ga=1$, and
the expected queue length $N=E(n)=1+E(|I|)$.
Once the stationary value of $\ga>0$ is reached, the probability distribution of $n$ is
given by
\be \label{5.4}
P_k=P(n=k)=P(|I|=k-\ga)
\ee
giving the obvious result that $k\ge\ga$. The explicit expression of the $P_k$ depends
therefore from the distribution of the $|I|$'s, and hence from the details of $f^{(\gs)}_\xi(t)$. 
\noindent
This solves completely the stationary problem in the $\gr=1$ case.
For application to the air traffic, however, it could be also interesting to study 
some non stationary features of the system: in particular we want to compute the probability
to pass from some negative value of $\ga$ to the following value $\ga+1$. 
These quantities are interesting in this $\gr=1$ case because if the probability to reach the state 
$n=0$ for a given $\ga\le 0$ is much smaller
that the inverse of the number of operation in a single day of traffic, it is
very likely that the system remains on states $n>0$. These probability to jump
from a definite value of $\ga$ to the following one are important also in the
description of the $\gr<1$ case, as it will be explained below.

\noindent
Hence suppose
that at time $j$ the system is in the state $n(j)=0$, with a given value of $\ga<0$.
Call $t(\ga)$ the quantity such that $n(j+i)>0$ for all $0<i<t(\ga)$, and
$n(j+t(\ga))=0$. $t(\ga)$ is therefore the length of the busy period with
starting value $\ga$. We are interested to the quantities
$T(\ga)=E(t(\ga))$.
By the definition of $\ga$ we have that $|I_j|=-\ga+1$ and that
the instant $j+t(\ga)$ is the first instant
after $j$ in which $|I_{j+t(\ga)}|=-\ga$, having $|I_{j+i}|>-\ga$ for all $0<i<t(\ga)$.
To compute $T(\ga)$ we should evaluate the probability
$P(|I_{j+i}|=-\ga{\big|}|I_j|=-\ga+1)$. This probability are however hard
to compute due to the conditioning.
Here we introduce our approximation: we will measure 
$T(\ga)$ in terms of 
\be \label{5.5}
T(\ga)\approx \frac{1}{P(|I|=-\ga)}
\ee
i.e. we neglect the conditioning. This approximation is reasonable for $\ga$
such that $P(|I|=-\ga)\ll\frac{1}{2L}$: 
in these cases we have to expect that the probability to have 
$P(|I_{j+i}|=-\ga{\big|}|I_j|=-\ga+1)$ for $i<2L$ is very small, and 
since $I_{j+i}$ is independent of $I_j$ for the 
greater values of $i$, that gives the bigger contribution to
$T(\ga)$, we have that the conditioning is almost ineffective.
On the other side, for $\ga$
such that $P(|I|=-\ga)\ge\frac{1}{2L}$ we have to expect a gross
underestimate of $P(|I_{j+i}|=-\ga{\big|}|I_j|=-\ga+1)$, and therefore
a gross overestimate of $T(\ga)$. We will return on this point later.

\noindent
We want now to compute explicitly $P(|I|=-\ga)$.
We will write general formulas, valid for any density $f^{(\gs)}_\xi(t)$,
and we will also consider a concrete probability distribution for the delays $\xi$,
namely the case of $f^{(\gs)}_\xi(t)$ uniform in $[-L,L]$, in which many computations
may be carried out explicitly.

\noindent
By straightforward computations one can see that 
\be \label{5.6}
P(|I|=0)= \prod_{i=-L+1}^{L-1}(1-F_{\xi}(i))=\frac{(2L)!}{(2L)^{2L}}\approx e^{-2L}\sqrt{4\pi L}
\ee
where the last approximation is valid for uniform $\xi$'s, using Stirling formula,
and
$$
P(|I|=k)= P(|I|=0)\sum_{-L+1\le i_1<i_2<...<i_k\le L-1}\frac{F_{\xi}(i_1)}{1-F_{\xi}(i_1)}...
\frac{F_{\xi}(i_k)}{1-F_{\xi}(i_k)}=
$$
\be \label{5.7}
=P(|I|=0)\sum_{-L+1\le i_1<i_2<...<i_k\le L-1}
\frac{L-i_1}{L+i_1}...\frac{L-i_k}{L+i_k}
\ee
where $F_{\xi}(t)$ is the probability distribution of the $\xi$'s, and the last
equality is again valid for uniform distribution.

\noindent
It is worthy to observe that \eqref{5.7} may be interpreted as the canonical partition
function of a Fermi system with $2L$ energy level and $k$ particles, where the $i$-th level has energy
$\log(F_{\xi}(i))-\log(1-F_{\xi}(i))$. With this respect many computational techniques
may be used in order to compute the probabilities $P(|I|=k)$. Note that, in the approximation 
\eqref{5.5}, once we are able to compute the quantities $P(|I|=k)$ we know also the
expected values $T(\ga)$.

\noindent
Let us list here a couple of possible way to evaluate $P(|I|=k)$ using the fact
that, since it is possible to interpret it as a well known object in statistical mechanics,
one can use computational results that are classical in that framework.
The number of energy level, as mentioned above, is $2L$. In real traffic
context one should expect that this value is of the order 20 or 30. One of the
available approximation of the quantity $P(|I|=k)$, i.e. the so called equivalence with
the \textit{grand canonical ensemble}, uses a method that is roughly speaking the 
Lagrange multipliers method, giving very good approximations for $2L$ large
(see e.g. \cite[chapter 5, section 53]{Landau}).
Since in our case $2L$ is not large enough to ensure the goodness of the
approximation, it is much better to use an exact expression for $P(|I|=k)$, due
to Ginibre. For completeness, and for the fact that it is quoted in a very
implicit sense in \cite{ginibre}, we give the proof of this formula.

\noindent
Calling $w_i=\frac{F_{\xi}(i)}{1-F_{\xi}(i)}$, one can prove the following equality
\be \label{canonical}
P(|I|=k)=\sum_{l=0}^k\sum_{1\le j_1\le...\le j_l\atop\sum_mj_m=k}C(j_1,...,j_l)\prod_{m=1}^l
\sum_i(w_i)^{j_m}
\ee
with
\be
C(j_1,...,j_l)=P(|I|=0)\frac{(-1)^{k-l}}{j_1.....j_lm_1!...m_k!}
\ee
where $m_i$ is the number of $j$'s equal to $i$.
To prove \eqref{canonical} we observe that
$$
\left.P(|I|=k)=P(|I|=0)\frac{1}{k!}\frac{d^k}{dt^k}\prod_i(1+tw_i)\right|_{t=0}
$$
The quantity $\prod_i(1+tw_i)$ can be expanded in series as follows
$$
\left.\frac{1}{k!}\frac{d^k}{dt^k}\prod_i(1+tw_i)\right|_{t=0}=\left.\frac{1}{k!}\frac{d^k}{dt^k}
e^{\sum_i\log(1+tw_i)}\right|_{t=0}=\left.\frac{1}{k!}\frac{d^k}{dt^k}
e^{\sum_i\sum_{j=1}^k(-1)^{j-1}\frac{(tw_i)^j}{j}}\right|_{t=0}=
$$
$$
=\left.\frac{1}{k!}\frac{d^k}{dt^k}
e^{\sum_{j=1}^k(-1)^{j-1}\frac{t^j}{j}\sum_i(w_i)^j}\right|_{t=0}=\left.\frac{1}{k!}\frac{d^k}{dt^k}
\sum_{l=1}^k\frac{(\sum_{j=1}^k(-1)^{j-1}\frac{t^j}{j}\sum_i(w_i)^j)^l}{l!}\right|_{t=0}=
$$
$$
=\sum_{l=0}^k\frac{(-1)^{k-l}}{l!}\sum_{ j_1,..., j_l\atop\sum_mj_m=k}\prod_{m=1}^l
\sum_i\frac{(w_i)^{j_m}}{j_m}=\sum_{l=0}^k\frac{(-1)^{k-l}}{l!}
\sum_{1\le j_1\le...\le j_l\atop\sum_mj_m=k}\prod_{m=1}^l
\sum_i\frac{(w_i)^{j_m}}{j_m}\frac{l!}{m_1!...m_k!}
$$
which is \eqref{canonical}.

\noindent
We conclude then the discussion of the $\gr=1$ case observing that in
a concrete framework of air traffic, if we want to avoid to have lost
slot but we want to keep the queue as short as possible we have to
choose initial condition in such a way that $\ga$ is the smaller
possible value such that $T(\ga)>D$, where $D$ is the number of
operations in a day. This value of $\ga$ gives the corresponding value
of the length of the queue using \eqref{5.4.1}. 

\noindent
A simple observation allows us to give an estimate of the average
length of the queue 
also when $\gr<1$. Let us suppose that we impose the condition $\gr<1$ keeping
the time between two expected arrivals equal to the service time, but assuming
that the arrivals are described by PSRA process with random deletion
(see \eqref{2.12.1}), with 
probability of deletion equal to $1-\gr$.  
It is easy to realize that this corresponds to say that
the value of $\ga$ has a probability $1-\gr$ to decrease by one.
Hence we have this picture of our queueing system: the queue is described
by a superposition of a slow varying process, the process that describes
the value of $\ga$, and a fast varying process, the one describing the $n$
for fixed $\ga$. If we are able to compute the distribution
probabilities of the 
values of $\ga$, we can evaluate the expected length of the queue (and even
its distribution) by \eqref{5.4.1}, weighted with the probabilities of
the various values of $\ga$.   

\noindent
In the unconditioned approximation \eqref{5.5}, the computation of the
stationary probabilities $\pi_\ga$ of $\ga$
is a standard task of the theory of the birth-and-death processes: the
evolution of $\ga$ is a discrete time birth-and-death process, with
transition probabilities  
$$P_{\ga,\ga'}=\begin{cases}
1-\gr\equiv \mu_{\ga} &\mbox{if } \ga'=\ga-1\\
P(|I|=-\ga)\equiv\gl_{\ga}&\mbox{if } \ga'=\ga+1\\
1-\gl_{\ga}-\mu_{\ga}&\mbox{if } \ga'=\ga\\
0&\mbox{otherwise}\end{cases}$$
and boundary conditions $\mu_{-L+1}=\gl_0=0$.  
We get the following linear system 
\begin{align*}
\gp_{-L+1}&=\gp_{-L+1}(1-\gl_{-L+1})+\gp_{-L+2}\gm_{-L+2}\\
\gp_i &=\gp_{i-1}\gl_{i-1}+\gp_{i+1}\gm_{i+1}+\gp_{i}(1-\gl_i-\gm_i)
\qquad -L+1 < i < 0 \\
\gp_0&=\gp_{-1}\gl_{-1}+\gp_{0}(1-\gm_0)
\end{align*}
whose solution is 
$$\gp_i=\gp_{-L+1}\prod_{k=-L+2}^{i}\frac{\gl_{k-1}}{\gm_{k}}
$$
The stationary distribution $\gp$ is defined by the normalization condition $\sum_i\gp_i=1$, then 
\be
\gp_{-L+1}=\frac{1}{1+\sum_{n=-L+2}^0
\prod_{k=-L+2}^{i}\frac{\gl_{k-1}}{\gm_{k}}
}
\ee 
This approximation is good
for $1-\gr$ sufficiently small,
because the probability to increase $\ga={-L+1}$ is much bigger than
the probability to decrease it, and at the same time the
unconditioned transition
probabilities to increase  $\ga$ when $\ga>{-L+1}$ 
are a good approximation of the actual transition probabilities.

\noindent
In the following figure we show the value of the expected
length of the queue obtained by the formula
\be
N=\sum_\ga\pi_\ga\left(\ga+E_\ga(|I|)\right)
\ee
Note that $E_\ga(|I|)$ is $\ga$-dependent, because in its
computation we neglect the terms with $|I|<-\ga$, since
they do not contribute to the evolution of the process with
that value of $\ga$. As it can be seen from the figure, the
estimate of the average length of the queue is extremely 
near to the simulations,
also for highly congested systems. In the figure we have
shown for completeness also the (wrong, for high $\gr$) values of the length
of the queue computed by means of formula \eqref{valore}, which neglects
the autocorrelations.

\begin{figure}[ht] 
\begin{center}
\includegraphics[width=.80\textwidth,totalheight=7cm]{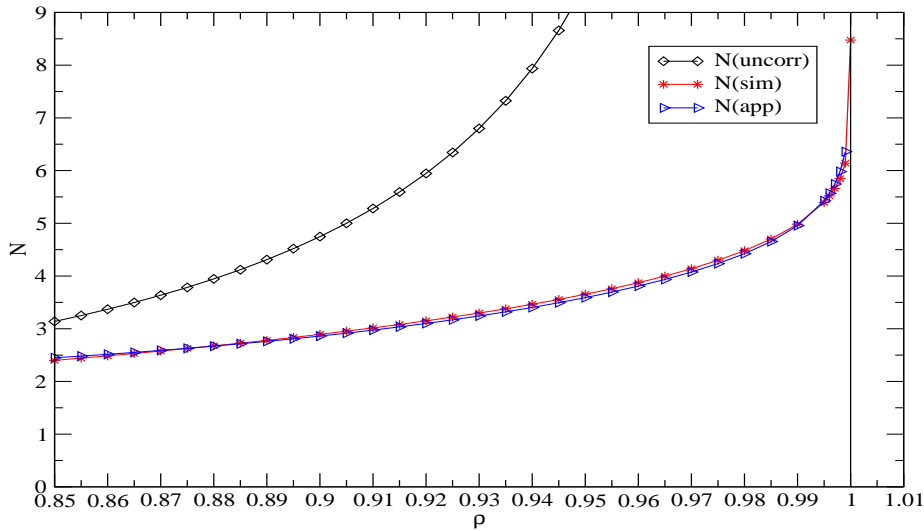}
\end{center}
\caption{The length of the queue for highly congested systems,
computed by means of numerical simulations (red line) and our 
analytical approximation (blue line). It can be seen that the
uncorrelated approximation (black line) obtained by formula \eqref{valore}
gives for these values of $\gr$ a gross overestimate. The simulations are run for a time sufficiently
long to have fluctuations on the result negligible in the scale of the figure.}
\label{lung}
\end{figure}
\noindent

\pagebreak
\section{Conclusions and open problems}

The main aim of this work is to study a stochastic process close to the Poisson process, 
but more suitable to describe the arrivals to a queueing systems when such arrivals
are scheduled in advance, and some randomness is added to the
schedule. We looked into this problem as an attempt to describe the
congestion in air traffic systems, but the same 
construction can be used in different contexts. 

We found analytical results, in particular we showed that our
process can be indistinguishable from a Poisson process if one wants
to study the distribution either of the number of arrivals or of the interarrival 
times in a time slot shorter
than the standard deviation of the randomness imposed to the scheduled
arrivals.

However we have shown that from the point of view of the resulting congestion,
due to the autocorrelation of this stochastic process, the queueing properties
of this model are quite different from the analogous problem with Poisson
arrivals. 
Interesting connection with the statistical mechanics emerged in the
analytical solution of the problem. We proposed some approximation
in our computations, but the results we obtained are in very good
agreement with numerical simulations.
An important question is the discussion of the accuracy of this
description with respect to actual air traffic data. We have
with this respect some preliminary results showing that
the description of the distribution of the length of the queue
using the PSRA as arrival process is much more accurate
than the description assuming Poisson process, that is
well known to be unfit.
We hope that this study, that has to be fully developed in its
computational aspects, may shed some light in various unclear aspects of
the air traffic modeling.

\section*{Acknowledgments}

S.N. is supported by Istituto Nazionale di
Alta Matematica ``Francesco Severi''. G.G. thanks the Dipartimento di
Matematica Universit\`a di Roma ``Tor
Vergata'' for its hospitality. We want to thank Domenico Marinucci and Giovanni
Peccati for discussions and encouragements. Errico Presutti gave us interesting
hints when we started to work on this problem. We thank Antonio Iovanella
for his help on the simulation. The Performance Review Unit of
Eurocontrol, and in particular Francesco Preti and Philippe Enaud suggested us
to investigate this problem.

\end{document}